\documentclass[11pt]{amsart} 

\usepackage{amsmath,amsfonts,amssymb,amsthm,amscd}
\usepackage{stmaryrd}

\usepackage{esint}

\usepackage{hyperref}
\usepackage{dsfont}

\usepackage{color}

\usepackage{bm}
\usepackage{amsbsy}
\usepackage{amssymb}
\usepackage[centertags]{amsmath}
\usepackage{amscd, amsthm}
\usepackage{stmaryrd}
\usepackage{enumitem}

\usepackage{xspace}

\usepackage{hyperref}

\newtheorem{prop}{Proposition}[section]

\newtheorem{thm}[prop]{Theorem}
\newtheorem{lemma}[prop]{Lemma}
\newtheorem{corollario}[prop]{Corollary}

\newtheorem{osserva}[prop]{Remark}

\theoremstyle{remark}

\def\be{\begin{equation}}
\def\ee{\end{equation}}

\newcommand{\B}{\mathbb B}

\newcommand{\R}{\mathbb R}
\newcommand{\Q}{\mathbb Q}

\def\dotminus{\mathbin{\ooalign{\hss\raise1ex\hbox{.}\hss\cr
  \mathsurround=0pt$-$}}}
\begin{document}
\title[Commutative unital rings elementary equivalent to ...]{Commutative unital rings elementarily equivalent to prescribed product rings } 
\author{Paola D'Aquino}
\address{Dipartimento di  Matematica e Fisica, Universit\`a della Campania L.Vanvitelli, 
Viale Lincoln 5, 81100 Caserta, Italy. 
E-mail: paola.daquino@unicampania.it}
\author{Angus J. Macintyre}
\address{School of Mathematics, University of Edinburgh, EH9 3FD Edinburgh, U.K. E-mail: a.macintyre@qmul.ac.uk}

\begin{abstract}
The classical work of Feferman-Vaught \cite{FV1959} gives a powerful, constructive analysis of
definability in (generalized) product structures, and certain associated enriched Boolean structures.
Here, by closely related methods, but in the special setting of commutative unital rings, we obtain
a kind of converse allowing us to determine in interesting cases, when a commutative unital R is
elementarily equivalent to a “nontrivial'' product of a family of commutative unital rings $R_i$. We use this in the model-theoretic analysis of residue rings of models of Peano Arithmetic. 
\end{abstract}

\maketitle
\section{Introduction} 
In several recent investigations (\cite{DM-ad2}, \cite{PAresrings}) the following question was considered:

When is a ring $R$ elementarily equivalent (in the standard first-order language $\mathcal L_{rings}$ for unital rings) to some {\it interesting} product $\prod_{i\in I}R_i$, where $R_i$ are commutative unital rings?

\medskip

We do not define {\it interesting}  here, but it is clearly not interesting to take $I=\{ 0\}$ and $R_0=R$. 

The question makes sense too for restricted products such as ad\`ele rings, for example $\mathbb A_{\Q}$, a restricted product of the rings $\Q_p$ ($p$ prime) and $\R$.

The two cases that stimulated our work are 
\begin{enumerate}[label=\roman*)]

\item
 rings $\mathcal M/n\mathcal M$, where $\mathcal M$ is the ring associated to a model of first-order Peano Arithmetic $PA$, in \cite{PAresrings}. The hardest case is when $n$ is nonstandard and divisible by infinitely many primes;
 
 \item
 classical ad\`ele rings $\mathbb A_{K}$, $K$ a number field, in \cite{DM-ad2}.

\end{enumerate}
The ring $\mathcal M/n\mathcal M$ is not a classical product of the rings $\mathcal M/q\mathcal M$, where $q$ runs through the maximal prime powers dividing $n$, but is an {\it internal product} (in the sense of nonstandard analysis). It turns out that $\mathcal M/n\mathcal M$ is elementarily equivalent to the full external product of the $\mathcal M/q\mathcal M$'s. This is rather surprising, and depends on refinements of the great work  of Feferman and Vaught from 1959, see \cite{FV1959}.

Our work on case ii) contributes  greatly to the understanding of the model theory of rings $\mathcal M/n\mathcal M$, a topic previously neglected (except when $n$ is prime, in \cite{MacResField}), but revived after   Zilber's \cite{BorisComm}.

\subsection{}

In this paper we refine a specially important case of the classical Feferman-Vaught Theorem (\cite[Theorem 3.1]{FV1959})
on products of first-order structures (some results of \cite{FV1959} were anticipated by Mostowski in \cite{Mos}). The
general case concerned a first-order language $\mathcal L$, and a family $\{R_i \}_{i \in I}$ of $\mathcal L$-structures, out of
which we form the classical $\mathcal L$-structure $R$ which is the full product of the $R_i$, and the objective
was to analyze definability (and decidability) in $R$ in terms of definability in the “stalks” $R_i$, and
in terms of the atomic Boolean algebra $\mathcal P(I)$. The method is amazingly general and
uniform, but has the defect that in general in \cite{FV1959} one has to appeal to various auxiliary Boolean
structures whose definability theory may be complicated. But even in the general setting
Weispfenning used the method (and extensions to generalized products) to obtain a wide variety of
important algebraic applications (\cite{Weispf}).
In this paper we restrict to the case when $\mathcal L_{rings}$ is the usual first-order language for rings, and the
 $R_i$ are unital, commutative rings, which are {\bf connected}  in the sense that $0$ and $1$
are the only idempotents and $0\not= 1$. We want a criterion (and in fact we get a necessary and
sufficient condition) for when a commutative unital ring $R$ is elementarily equivalent to some
product of connected rings $R_i$ as above.
What we get is a weak converse to  Theorem 3.1 of \cite{FV1959}. That theorem is
of much greater generality than anything we need here, with many notational complications. We
work entirely in the language $\mathcal L_{rings}$ of ring theory, coding  the auxiliary Boolean formalism of \cite{FV1959} in terms
of ring idempotents, and in particular minimal idempotents (the coding is always done in the ring $R$
where we are working). The analysis of  Theorem 3.1 of \cite{FV1959}  in our ring theoretic setting leads to the isolation of a set AXIOMS in the ring language, and 
a proof that products of connected rings satisfy the axioms we isolate. 
Our main result Theorem \ref{FVunital} shows that the models of  AXIOMS  are exactly the rings which are elementarily equivalent to a product of connected unital rings. 
We have found no such “converse'' in \cite{FV1959}.

\

\noindent
{\it Examples of connected rings} 1. Every integral domain is  connected. 

\smallskip

\noindent 2. Any commutative local ring is connected, as we now show. If $S$ is local and $x\in S$ satisfies $x^2=x$ then either $x\in \mu$ or $x-1\in \mu$, where $\mu$ is the unique maximal ideal of $S$.  If $x\not\in \mu$ then $x$ is a unit, and so $0=ux(x-1) = x-1$ where $u$ is the inverse of $x$. If $x-1\not\in \mu$ then $x-1$ is a unit, and hence $0=x(x-1)u$ where $u$ is the inverse of $x-1$. So, either $x=1$ or $x=0$. It cannot happen that $x,x-1\in \mu$ otherwise $1\in \mu$, which is a contradiction.

\subsection{}
The Feferman-Vaught analysis of products and generalized products of first order $\mathcal L$-structures in \cite{FV1959} is done in terms of  $\mathcal L$ and also the language of Boolean algebras and expansions of it. 

In the case of products of connected commutative unital rings $R=\prod_{i\in I}R_i$  the hypothesis that the $R_i$ are connected is used to make the following identification  of the Boolean algebra of the power set of $I$ in ring theoretic terms inside $R$.  The elements of $R$ are functions $f$ defined on $I$ such that $ f(i)\in R_i$, and the  Boolean algebra $\B$ associated to $R$ is the set of functions in $R$  which take values in $\{0,1\}$, i.e. the set of idempotents of $R$ (here the hypothesis of connectness of $R_i$ is crucial), with a natural Boolean structure defined on it (see the next section).  The individual element $i\in I$ is identified with the function $\chi_i$  with $\chi_i(i)=1$, and $\chi_i(j)=0$ for $j\in I$ and $j\not=i$,  and these correspond to the atoms of $\B$. 

In the case of a general commutative unital ring $R$ (not a product of connected commutative unital rings) we will use the same  notation $\B$   to denote the Boolean algebra of idempotent elements of $R$, and $e$ will be an atom of $\B$. Note that for $R=\prod_{i\in I}R_i$ ($R_i$ connected commutative unital rings) the Boolean algebra $\B$ is complete while in the general case $\B$ may not be complete. Inspired by the Feferman-Vaught arguments in the proof of Theorem 3.1  in \cite{FV1959} we will identify axioms, mainly in terms of $\B$, for $R$ in order to have $R$ elementarily equivalent to a product of connected commutative unital rings which will be determined in terms of the atoms of  $\B$.



We will find these properties, and prove a theorem applicable to both cases (i) and (ii) of page 1.  

Note that when $I$ is infinite and each $R_i$ is nontrivial then $\prod_{i\in I}R_i$ is uncountable.

\section{Rings with connected stalks: a Feferman-Vaught theory}
\subsection{Basic idea}
We will work in the languages of unital rings $\mathcal L_{rings}$ which contains $+, -, \cdot, 0,1$,  and  in the language of Boolean algebras $\mathcal L_{B}$ which contains $\pmb{\wedge}, \pmb{\vee}, \pmb{\neg}, 0,1 $.  Let $R$ be any commutative   unital ring, $\B$ the Boolean algebra of idempotents of $R$, interpretable in $R$ via

\medskip

$
\begin{array}{ccc}
e\pmb{\wedge} f &   =    & ef   \\
e\pmb{\vee} f   &    =   &  1-(1-e)(1-f)   \\
              &    =     &    e+f-ef \\
\pmb{\neg} e &   =    &   1-e   \\
0    &    =    &    0  \\
1    &    =    &    1    
\end{array}
$

\medskip

One may note that on the set $\B$ we have the equivalent Boolean ring structure $\B_{ring}$ (generally, not a subring of $R$) with

\medskip

$
\begin{array}{ccc}
x\oplus y &   =    & (x{\pmb{\wedge}} \pmb{\neg} y)\pmb{\vee} (\pmb{\neg}  x\pmb{\wedge} y)  \\
x\odot y   &    =   &  x \pmb{\wedge} y  \ ( =  x\cdot y)  \\
0    &    =    &    0  \\
1    &    =    &    1    
\end{array}
$

\medskip

On $\B$ we have the usual definable partial order $\leq $. Of special importance to us are the atoms of $\B$ (equivalently, the minimal idempotents $\not= 0, 1$).

\begin{osserva}
\label{onesort}
 {\rm The set of idempotents of a ring $R$ is quantifier-free definable in the ring language and via the above interpretation of $\pmb{\wedge}, \pmb{\vee}, \pmb{\neg}$ the Boolean structure on $\B$ is quantifier-free  in terms of $+, -, \cdot, 0,1$ in $R$. In this way any formula in $\mathcal L_{B}$ can be  effectively {\it translated} into a formula of $\mathcal L_{rings}$. We will refer to these formulas as {\it B-formulas}. In this way we are going to work in a one sorted language as opposed to the language in \cite{FV1959}  which is a two sorted language, one for the family of structures indexed by $I$, and the other of the Boolean algebra. }
\end{osserva}

\begin{lemma}
\label{localized}
For any idempotent element $e$ in $\mathbb B$ we have natural isomorphisms 
$$R/(1-e)R \cong eR \cong R_e, $$
where $R_e$ is the localization of $R$ at $\{ e^n: n\geq 0\}$ $(=\{ 1, e\})$. 
\end{lemma}

\proof
The ring $R$ is the direct sum of the ideals $eR$ and $(1-e)R$, giving the first isomorphism. The kernel of the map $R\rightarrow  R_e$ consists of the $x$ so that $x
e=0$, and this is the same as $x=x
e+ x
(1-e)=x
(1-e)$, thus the same as $x\in (1-e)R$. The ring $eR$ is unital with unit $e$. Using the universal property of localization we can explicitly define the isomorphism between $R_e$ and the ring $eR$ as $\frac{r}{s} \mapsto er$, where $s=1,e$. 
\endproof

\noindent 
{\bf Notation. } An element of $R_e$  will be denoted by $f_e$ for some $f\in R$, and by Lemma \ref{localized} it corresponds to $ef$ in $eR$.

By an atom in $R$ we will always mean an atom in the Boolean algebra $\mathbb B$ of idempotents of the ring $R$, and the order $<$ is the order on $\mathbb B$. 

\begin{corollario}
A nonzero idempotent $e$ of $R$ is an atom if and only if $R_e$ is connected.
\end{corollario}

\proof
\noindent ($\Rightarrow$) If $e$ is an atom in $R$ then $R_e$ is the  ring with two elements which is clearly connected. 

\noindent ($\Leftarrow$) Let $e$ be a nonzero idempotent of $R$, and assume $R_e$ is connected. We show that $e$ is an atom in $R$. Let $f\in R$ be an idempotent, and $0<f<e$. Then $fe=f$ (by the Boolean operations on $\mathbb B$), and so $f$ is not in the kernel of the localizing map at $e$ from $R$ to $R_e$ since $f\not=0$. By connectedness of $R_e$, $f=e$ in $R_e$, and so $e(f-e)=0$ in $R$, which implies $ef=e$. Hence, $f=e$ in $R$, a contradiction. 
\endproof


We have defined the localizations of $R$ at any idempotent, but the interesting case is when $e$ is an atom. Only in this case we have that $R_e$ is connected, and this is a fundamental property of the stalk $R_e$ of $R$.  We will identify axioms which characterize commutative unital rings which are elementarily equivalent to the product of their stalks at atoms. 

\smallskip

\medskip

\noindent
{\bf Axiom 1.} The Boolean algebra $\B$ is atomic.

\medskip

This clearly holds in any product of commutative unital rings, and so in particular it holds in products of commutative connected unital rings. Notice that there are commutative rings where Axiom 1 fails, in particular take $R$ the Boolean ring associated to an atomless Boolean algebra. 

\smallskip

The following property is well known, see e.g. \cite{GH}. We will often use it in the proofs of our results.

\medskip




\begin{lemma}
\label{vee}
If $f\in \B$, $f\not=0$ then $$f=\bigvee_{\substack{e\leq f \\e \mbox{ {\tiny atom}}}}e.$$

\end{lemma}

The proof of the lemma follows easily from the following equivalences that are proved on page 118 of \cite{GH}
\begin{enumerate}
\item
a Boolean algebra $B$ is atomic;

\item
every element $b\in B$ is the supremum of the set of atoms $a\leq b$;
\item
$1\in B$ is the supremum of all atoms.
\end{enumerate}

\subsection{Completeness}
In a product   $R=\prod_{i\in I}R_i$ where $R_i$'s are connected  commutative unital rings, the Boolean algebra $\B$ is complete. This means that every non empty subset of $\B$ has a supremum (and hence also an infimum) with respect to the usual order on the Boolean algebra. 

This is not a first order condition, and we need to find a replacement. We use a variant, adapted to ring theory. 


For any $f_0, \ldots, f_n\in R$, and any  atom $e\in \mathbb B$ we will use the notation $\overline{f}$ for the tuple $f_0, \ldots, f_n$, and     $(\overline{f})_e$   for  the tuple $(f_{0})_e, \ldots, (f_{n})_e$ of elements of the stalk $R_e$.


\noindent
{\bf Axiom 2.} For all formulas $\theta (x_0, \ldots, x_n)  $ in $\mathcal L_{rings}$,  and for all  $f_0, \ldots, f_n\in R$ there is a unique element $\llbracket   \theta (f_0, \ldots, f_n)   \rrbracket$ in $\mathbb B$ such that 
$$R_e \models  \theta ((f_{0})_e, \ldots, (f_{n})_e) \mbox{ \hspace{.1in} iff  \hspace{.1in}} e \leq \llbracket   \theta (f_0, \ldots, f_n)   \rrbracket$$ for all atoms $e$ in $\mathbb B$.

\medskip
The elements $\llbracket   \theta (f_0, \ldots, f_n)   \rrbracket$ occur in \cite{FV1959} with a different notation (see Definition 2.1 on page 63 of  \cite{FV1959}). The $\llbracket  \ \cdot \  \rrbracket$ notation comes from Boolean-valued model theory \cite{Bell}, not yet prominent when \cite{FV1959} was written.

 Notice that Axiom 2 is a scheme of axioms, one for each formula $\theta (x_0, \ldots, x_n)$. This is our substitute for completeness. 
 
 Axiom 2  is clearly true in products of unital connected commutative rings, since $\llbracket   \theta (f_0, \ldots, f_n)   \rrbracket$ is  $\bigvee \{ i: R_i \models \theta (f_{0}(i), \ldots, f_{n}(i))\}$.


 
 \medskip

\subsection{Boolean values}

We need to check preservation properties of $\llbracket   \theta (f_0, \ldots, f_n)   \rrbracket$. We will use Axiom 2  in the proofs.
\begin{lemma}
\label{conjunction}
Let $\theta_1, \theta_2$ formulas, and $f_0, \ldots, f_n \in R$. Then $$\llbracket   (\theta_1\wedge  \theta_2)(\overline{f} )   \rrbracket = \llbracket   \theta_1(\overline{f} )\rrbracket \wedge  \llbracket \theta_2(\overline{f} )   \rrbracket,$$ where $\overline{f}=f_0, \ldots, f_n$.

\end{lemma}

\proof
Let $e$ be an atom. The following equivalences hold:

\smallskip

\noindent $e\leq \llbracket   (\theta_1\wedge  \theta_2)(\overline{f} )   \rrbracket$ iff $R_e\models (\theta_1\wedge  \theta_2)((\overline{f})_e )$ iff   $R_e\models \theta_1((\overline{f})_e )$  and $R_e~\models  ~\theta_2~((\overline{f})_e~)$ iff $e\leq \llbracket   \theta_1 ((\overline{f})_e)   \rrbracket$ and $e\leq \llbracket   \theta_2 ((\overline{f})_e)   \rrbracket$ iff $e\leq \llbracket   \theta_1 ((\overline{f})_e)   \rrbracket  \wedge  \llbracket   \theta_2 ((\overline{f})_e)   \rrbracket$. 

The result now follows from Lemma \ref{vee}. \endproof

\begin{lemma}
\label{negation}
Let $\theta$ be a formula, and $f_0, \ldots, f_n \in R$. Then $$\llbracket   \neg \theta (\overline{f})   \rrbracket= \neg  \llbracket   \theta (\overline{f})   \rrbracket.$$ 

\end{lemma}

\proof
 Recall that for any atom $e$ and any $\alpha$ in a Boolean algebra either $e\leq \alpha $ or $e\leq \neg \alpha$, and not both inequalities hold. Let $e$ be an atom. The following equivalences hold:

\smallskip

\noindent
$e\leq   \llbracket   \neg \theta (\overline{f})   \rrbracket$ iff  $R_e\models \neg \theta ((\overline{f})_e )$ iff   $R_e\not\models \theta ((\overline{f})_e )$ iff 
$e \not\leq \llbracket   \neg \theta (\overline{f})   \rrbracket$ iff    

\noindent $e \leq \neg \llbracket   \theta (\overline{f})   \rrbracket$.  

 The result now follows from Lemma \ref{vee}. \endproof

\begin{lemma}
\label{disjunction}
Let $\theta_1, \theta_2$ be formulas, and $f_0, \ldots, f_n \in R$. Then $$\llbracket   (\theta_1\vee  \theta_2)(\overline{f} )   \rrbracket = \llbracket   \theta_1(\overline{f} )\rrbracket \vee  \llbracket \theta_2(\overline{f} )   \rrbracket,$$ where $\overline{f}=f_0, \ldots, f_n$.

\end{lemma}

\proof
This follows from the two preceding lemmas. \endproof

\medskip

\noindent
{\bf Axiom 3.} For all formulas $\theta (v_0, \ldots, v_n,w)$ in $\mathcal L_{rings}$, and for all $f_0, \ldots, f_n\in R$ there is $g\in R$ such that 
$$ \llbracket  \exists w \theta (\overline{f},w)     \rrbracket \leq \llbracket   \theta (\overline{f},g)     \rrbracket 
$$
where $\overline{f}=f_0, \ldots, f_n$.

\begin{lemma}
Axiom 3 is true in $\prod_{i\in I}R_i$ where $R_i$'s are connected  commutative unital rings.
\end{lemma}

\proof Let $h= \llbracket  \exists w \theta (\overline{f},w)  \rrbracket \in \mathbb B$ be defined as  
$$h(i)= \left\{ \begin{array}{ll}
1 & \mbox{ if } R_i \models \exists w\theta (f_{0}(i), \ldots, f_{n}(i), w) \\
0  & \mbox{ if } R_i \not\models \exists w\theta (f_{0}(i), \ldots, f_{n}(i),w)
\end{array} \right.$$

Now we define $g \in \prod R_i$ as follows 

$$g(i)= \left\{ \begin{array}{ll}

g_i & \mbox{ if } R_i \models \theta (f_{0}(i), \ldots, f_{n}(i), g_i) \\
0  & \mbox{ if } R_i \not\models \exists w\theta (f_{0}(i), \ldots, f_{n}(i),w)
\end{array} \right.$$

We use the axiom of choice in order to pick a $g_i\in R_i$  such that  $h(i)=1 $ and  $R_i \models \theta (f_{0}(i), \ldots, f_{n}(i),g_i) $.   If $R_i \not\models \exists w\theta (f_{0}(i), \ldots, f_{n}(i),w) $ then  $g(i)$ can be chosen arbitrarily, we opted to define it $0$. Let $t= \llbracket   \theta (\overline{f},g)     \rrbracket \in \mathbb B$ be defined as follows 
$$t(i)= \left\{ \begin{array}{ll}
1 & \mbox{ if } R_i \models \theta (f_{0}(i), \ldots, f_{n}(i), g(i)) \\
0  & \mbox{ if } R_i \not\models \theta (f_{0}(i), \ldots, f_{n}(i),g(i))
\end{array} \right.$$
It is immediate to show that  $h\leq t$. \endproof

\subsection{An axiom about atomic formulas}
We add a further axiom, clearly true in products, and henceforward assumed. 

\smallskip

\noindent
{\bf Axiom 4.} For any atomic formula $\theta (v_0, \ldots, v_n)$ in $\mathcal L_{rings}$,
$$ R \models \theta (\overline{f})   \mbox{\hspace{.1in}  iff \hspace{.1in} }   \mathbb B \models  \llbracket   \theta (\overline{f})   \rrbracket=1,  $$ where $\overline{f}=f_0, \ldots, f_n\in R. $

\medskip
Axiom 4  is clearly true in product of connected rings since the ring structure on $\prod_{i\in I}R_i$ is defined coordinatewise.

\subsection{Partitions and a Patching Axiom}



In order to sketch a proof of a useful generalization of \cite{FV1959} to our more restricted setting of rings $R$ satisfying the axioms listed above, we need to review several notions of partitions used in \cite{FV1959}. 

\smallskip
 In a Boolean algebra $B$ a partition is a finite sequence $Y_0,\ldots , Y_m$ of elements of $B$ such that $1=Y_0\vee \ldots \vee Y_m$, and $Y_i\wedge Y_j=0$ if $i\not=j$. (We do not insist that each $Y_i\not=0$, but we do insist that the sequence is finite.) This notion can be easily formalized in $\mathcal L_B$ (and in our context also in $\mathcal L_{rings}$) by the formula
  $$Part_{m+1}(y_0, \ldots, y_m):= (\bigvee_{0\leq j\leq m} y_j)=1 \wedge \bigwedge_{0\leq i<j\leq m}(y_i\wedge y_j=0).$$
  In our context via the interpretation of the Boolean algebra of idempotents in the ring $R$ we can give a meaning to the notation $R\models \Phi (  \llbracket   \theta_0     \rrbracket ,\ldots ,   \llbracket   \theta_m    \rrbracket )$, where  $\Phi$ is the $B$-formula of  $\mathcal L_{rings}$   corresponding to a formula  $\mathcal L_B$ (via the interpretation),  and $\theta_0, \ldots, \theta_m $ are  in $\mathcal L_{rings}$.
  Feferman and Vaught work in a two sorted language, one for the structure and one for the Boolean algebra,  while we have only a one sorted language without the Boolean sorts (see Remark \ref{onesort}). They introduce the following notions which occur in Theorem 3.1 in \cite{FV1959}. We recall them and adapt them if necessary to our context. 
  
\noindent 
1. An  {\it acceptable  sequence} is a $\zeta =\langle  \Phi, \theta_0, \ldots, \theta_m\rangle$ where  $\Phi $ is a $B$-formula in $\mathcal L_{rings}$ and $\theta_0, \ldots, \theta_m$ are $\mathcal L_{rings}$-formulas. From this one gets (after imposing trivial syntactic constraints) an $\mathcal L_{rings}$-formula $\Phi ( \llbracket \theta_0(\overline{v}) \rrbracket , \ldots,  \llbracket \theta_m(\overline{v})  \rrbracket ) $ which defines in any ring $R$ satisfying the axioms we have given so far a subset of $R^{|\overline{v}|}$. Note that $\zeta$ can be considered as a code (i.e. a natural number)  for $\Phi ( \llbracket \theta_0(\overline{v}) \rrbracket , \ldots,  \llbracket \theta_m(\overline{v})  \rrbracket ) $ by  associating a natural number to each  symbol of the language page 64  in  \cite{FV1959}. 

\noindent 
2. By a free variable in a sequence $\zeta $ we mean a variable $v$ that appears free in at least one $\theta_0, \ldots, \theta_m$. A {\it standard acceptable  sequence} is an acceptable  sequence $\zeta =\langle  \Phi, \theta_0, \ldots, \theta_m\rangle$ where the free variables are exactly $v_0, \ldots , v_k$ for some $k\geq 0$.

\noindent 
3. A sequence of $\mathcal L_{rings}$-formulas $\theta_0(\overline{v}), \ldots, \theta_m(\overline{v})$ is called a {\it partition sequence } if for any ring $R$ satisfying the axioms introduced so far and for any $\overline{f}$ elements of $R$ the sequence 
$$\llbracket \theta_0(\overline{f}) \rrbracket , \ldots,  \llbracket \theta_m(\overline{f})  \rrbracket $$     forms a partition of the Boolean algebra $\B$ of idempotents of $R$.  

\medskip

Feferman and Vaught  remark that in the language obtained  by adding all the $\zeta $ to the one sorted language of the structures $\mathcal M_i$ for $i\in I$ the product $\prod_{i\in I}\mathcal M_i$ has an elimination of quantifiers. 

We now give an idea on how to produce partition sequences following page 67 of  \cite{FV1959}. This will be used in the inductive step of Theorem  \ref{FVunital} when we have to consider existential formulas.

Given an acceptable sequence $\zeta' =\langle  \Phi', \theta_0', \ldots, \theta_{m'}'\rangle$ we can determine in an effective way another acceptable sequence $\zeta =\langle  \Phi, \theta_0, \ldots, \theta_m\rangle$ with the same free variables such that $\theta_0, \ldots, \theta_m$ form a partition sequence and for all models $R$ of our axioms 
\begin{equation} 
\label{Phi}
R\models \Phi ( \llbracket \theta_0(\overline{f}) \rrbracket , \ldots,  \llbracket \theta_m(\overline{f})  \rrbracket ) \mbox{  iff  } R\models \Phi' ( \llbracket \theta_0'(\overline{f}) \rrbracket , \ldots,  \llbracket \theta_{m'}'(\overline{f})  \rrbracket ) 
\end{equation}
  for all $\overline{f} \in R$. The acceptable sequence $\zeta$ is constructed as follows: 

\smallskip
\noindent 1. let $m=2^{m'+1}$ and $r_0, \ldots r_m$ be the list of subsets of $  m'+1=\{ 0, \ldots, m'\}$. 

For each $k\leq m$ define $$\theta_k := \bigwedge_{k\in r_k}\theta_k' \wedge  \bigwedge_{k\not\in r_k}\neg \theta_k';$$

\smallskip

\noindent
2. for each $ l \leq m'$ let $s_l:=\{ k: k\leq m \mbox{ and } l\in r_k\}$  and  $$\Phi (x_0, \ldots, x_m):=    \Phi' (\bigvee_{ k\in s_0}x_k, \ldots, \bigvee_{ k\in s_{m'}} x_k).   $$

It is easy to check that the $\theta_k$ form a partition, and that (\ref{Phi}) holds. Note that this argument is nothing more than the usual procedure of constructing the disjunctive normal form of a propositional logic formula, where the $\theta_k $  play the role of  propositional variables. 

We next have to introduce a new axiom which will use in the inductive step of the proof of our main result, Theorem \ref{FVunital}, when dealing with the existential case. This axiom is necessary since we are working in a more general setting than product of structures, and in our context  we cannot prove the properties in equation (15)  on page 68 of \cite{FV1959}. This will be the content of our fifth and  final  axiom.

 Axiom 5, 
 gives sufficient conditions for a ubiquitous “Patching Argument”, and it would surely be worthwhile to give a more sheaf-theoretic, or topos-theoretic, formulation.

\smallskip

Let $\phi (v_0,….v_m)$ be any formula  in  the language $\mathcal L_B$ of Boolean algebras with free variables $v_0,\ldots,   v_m$. To  $\phi (v_0,….v_m)$   we associate a formula $\phi^*$   with free variables $v_0,\ldots, v_m$  in $\mathcal L_B$ expressing:

there exists a partition $W_0, \ldots,   W_m$  of $B$ such that $W_j \leq v_j$ for each $j=0,\ldots, m$, and  $\phi (W_0, \ldots, W_m) $ holds in $B$.

\medskip

\noindent
{\bf Axiom 5.} Let  $\phi (v_0,….v_m)$ be a formula in  the language $\mathcal L_B$ of Boolean algebras with free variables $v_0,\ldots,   v_m$, and $\phi^*(v_0,\ldots, v_m)$ as above. 

Let $\theta_0(x_0,\ldots ,x_k, x_{k+1}), \ldots, \theta_m(x_0,\ldots ,x_k, x_{k+1})$ be formulas in the ring language and assume that $$(\theta_0(x_0,\ldots ,x_k, x_{k+1}), \ldots, \theta_m(x_0,\ldots ,x_k, x_{k+1}))$$ is a partition. Then the following  are equivalent in $\mathbb B$  for $k$-tuples $\bar{f}\in R$ 

(1) there exists $g\in R$ such that $  \phi^* (\llbracket \theta_0(\bar{f},g)\rrbracket , \ldots,\llbracket \theta_m(\bar{f}, g)\rrbracket ) $ holds in $\mathbb B$;

(2) there is a partition $Y_0,\ldots , Y_m$ of $\mathbb B$  such that for each $j~=~0~,~\ldots~, m$, $Y_j \leq \llbracket  \exists x_{k+1}\theta_j (\bar{f}, x_{k+1})\rrbracket $  and 
 $ \phi (Y_0,\ldots , Y_m) $ holds in $\mathbb B$.

\begin{lemma}
\label{axiom5}
Let $R=\prod_{i\in I}R_i$, where  $R_i$ are connected unital rings and  $\B$ the boolean algebra of idempotents of $R$.
Then Axiom 5 holds in $R$.

\end{lemma}

\proof
We use the same notation as in  Axiom 5. Recall  that in this case $\mathbb B$ can be identified with $\mathcal P(I)$, and  a partition of $\mathbb B$ is a finite collection of disjoint subsets of $I$ whose union is $I$. Moreover,  $$\llbracket  \xi (\bar{f})\rrbracket =\bigvee \{ i\in I : R_i\models \xi(f_0(i), \ldots , f_k(i)\}$$ for any $\mathcal L_{rings}$-formula $\xi(\bar{x})$ and any $\bar{f}\in  \prod_{i\in I}R_i$.

\smallskip

2) $\Rightarrow$ 1) Let $Y_0,\ldots , Y_m\in \mathbb B $ be a partition of $\B$ such that for each $j~=~0~,~\ldots~, m$, $Y_j \leq \llbracket  \exists x_{k+1}\theta_j (\bar{f}, x_{k+1})\rrbracket $  and $\B \models \phi (Y_0,\ldots , Y_m) $. Notice that for each atom $i\in \mathbb B$ there is a unique $j\in \{1,\ldots, m\}$ such that  $i\leq Y_j$ (this follows from the hypothesis that $Y_0,\ldots , Y_m$ is a partition of $\B$).  Fix $j$,  for all atoms $i$ such that $i\leq Y_j\leq \llbracket  \exists x_{k+1}\theta_j (\bar{f}, x_{k+1})\rrbracket $, we have $R_i \models  \exists x_{k+1}\theta_j (f_{0}(i), \ldots, f_{n}(i), x_{k+1})$. Using the axiom of choice we are now going to define an element $g\in \prod_{i\in I}R_e$ as follows: 

for all atoms $i\leq Y_j$ we set $g(i) = g_i\in R_i$ for some $g_i\in R_i$ such that  $ R_i\models  \theta_j (f_{0}(i), \ldots, f_{n}(i), g_i)$. We do this for all $j=1, \ldots ,m$, and so $g$ is well defined on all $i\in I$. 



Clearly,  $Y_j\leq \llbracket   \theta_j (f_0, \ldots, f_k,g)   \rrbracket$ for $j~=~0~,~\ldots~, m$. Indeed, $Y_j= \llbracket   \theta_j (f_0, \ldots, f_k,g)   \rrbracket$. 
If $i\leq \llbracket   \theta_j (f_0, \ldots, f_k,g)   \rrbracket$ 
 and $i\not\leq Y_j$ then $i\leq Y_h$ for a unique $h\not=j$, and $i\leq Y_h\leq  \llbracket   \theta_h (f_0, \ldots, f_k,g)   \rrbracket$. Then  $\llbracket   \theta_j (f_0, \ldots, f_k,g)    \wedge    \theta_h (f_0, \ldots, f_k,g)   \rrbracket = \llbracket   \theta_j (f_0, \ldots, f_k,g)   \rrbracket \wedge  \llbracket   \theta_h (f_0, \ldots, f_k,g)   \rrbracket \not =0$ which is a contradiction since $(\theta_0(x_0,\ldots ,x_k, x_{k+1}), \ldots, \theta_m(x_0,\ldots ,x_k, x_{k+1}))$ is a partition.
Hence, $\B \models \phi^* (\llbracket   \theta_0 (\bar{f},g)   \rrbracket,\ldots , \llbracket   \theta_m (\bar{f},g)   \rrbracket) $. 

\smallskip

1) $\Rightarrow$ 2) Suppose there is $g\in R$  such that $$\B \models \phi^* (\llbracket \theta_0(\bar{f},g)\rrbracket , \ldots,\llbracket \theta_m(\bar{f}, g)\rrbracket ) .$$ So, there is a partition $w_0, \ldots , w_m$  of $\B$ and $w_j\leq \llbracket \theta_j(\bar{f}, g)\rrbracket$ for all $j~=~0~,~\ldots~, m$, and $\B\models \phi (w_0, \ldots ,w_k)$. Recall that $$\llbracket \theta_j(\bar{f},g) \rrbracket  \leq \llbracket \exists x_{k+1}\theta_j(\bar{f},x_{k+1})\rrbracket $$ for all $j=1, \ldots , m$. Then 2) easily follows. \endproof

\






The set of the five axioms will be denoted by $AXIOMS$. From $AXIOMS$ we will prove our main result, Theorem \ref{FVunital} below. It is inspired by the fundamental Theorem 3.1 of \cite{FV1959}, but it is considerably less general, and uses quite different notation.



\begin{thm}
\label{FVunital}
 For every $\mathcal L_{rings}$-formula $\theta (x_0,\ldots ,x_k)$ there is a partition $(\theta_0(x_0,\ldots ,x_k), \ldots ,\theta_m(x_0,\ldots ,x_k))$  of ring formulas, and a Boolean algebra formula $\psi (y_0,\ldots , y_m)$ so that for all $f_0, \ldots ,f_k\in R$, where $R$ is a ring satisfying $AXIOMS$,  and $\mathbb B$ is the Boolean algebra of idempotents of $R$
 \begin{equation} 
\label{equivMain}
R\models \theta (\overline{f}) \mbox{ \hspace{.2in} iff \hspace{.2in}} \mathbb B \models \psi (\llbracket  \theta_0 (\overline{f}) \rrbracket, \ldots ,\llbracket  \theta_m (\overline{f}) \rrbracket )
\end{equation} 
 where $\overline{f}=f_0, \ldots ,f_k.  $
\end{thm}

\proof
As in \cite{FV1959} the proof proceeds by induction on the complexity of the formula $\theta (x_0,\ldots ,x_k)$. Using the five axioms we have identified we can easily prove the cases of atomic formulas (use Axiom 4) and the Boolean connectives (use Axiom 2 and Lemmas \ref{conjunction} and \ref{negation}).  When $\theta (x_0,\ldots ,x_k)=\exists x_{k+1} \xi (x_0,\ldots ,x_k, x_{k+1})$ Axiom 5 guarantees that there are no contradictions in the patching argument in the inductive step. Let $$(\xi_0(x_0,\ldots ,x_k,x_{k+1}), \ldots ,\xi_m(x_0,\ldots ,x_k,x_{k+1}))$$ be the partition in $\mathcal L_{rings}$-formulas and $\psi (v_0,\ldots ,v_m)$ the formula in $\mathcal L_B$ associated by the inductive hypothesis to $\xi (x_0,\ldots ,x_k, x_{k+1})$.  Using Axiom 5  and the inductive hypothesis the following are equivalent
\begin{enumerate}
 \item $R\models \exists x_{k+1} \xi (f_0,\ldots ,f_k, x_{k+1})$
 
\item there exists $g\in R$ such that  $R\models \xi (f_0,\ldots ,f_k, g)$
 
\item $\mathbb B\models \psi  (\llbracket  \xi_0 (\overline{f},g) \rrbracket, \ldots ,\llbracket  \xi_m (\overline{f},g) \rrbracket )   $
\item $\mathbb B\models \psi'  (\llbracket  \exists x_{k+1}\xi_0 (\overline{f},x_{k+1}) \rrbracket, \ldots ,\llbracket  \exists x_{k+1}\xi_m (\overline{f},x_{k+1}) \rrbracket )   $ where 
$\psi'  (v_0, \ldots ,v_m) $ is the $\mathcal L_B$-formula $$ \exists w_0, \ldots, w_m (\bigwedge_{i=0}^m  w_i\leq v_i\wedge  (w_0, \ldots, w_m  \mbox{ is partition of $\mathbb B$}) \wedge \psi ( w_0, \ldots w_m)).$$
\end{enumerate} 
Notice that $\exists x_{k+1}\xi_0 (x_0,\ldots ,x_k,x_{k+1}), \ldots , \exists x_{k+1}\xi_m (x_0,\ldots ,x_k,x_{k+1})$ may not be a partition. As outlined on page 67 of \cite{FV1959} (and remarked also in this paper before introducing axiom 5) there is an effective procedure which associates a partition to  $\exists x_{k+1}\xi_0 (\overline{x},x_{k+1}), \ldots , \exists x_{k+1}\xi_m (\overline{x},x_{k+1})$ mantaining the same free variables.

 \endproof

 The statement of Theorem \ref{FVunital} is uniform for all models of our axioms, and this uniformity  is fundamental. The partition $$(\theta_0(x_0,\ldots ,x_k), \ldots ,\theta_m(x_0,\ldots ,x_k))$$  of ring formulas, and the Boolean algebra formula $\psi (y_0,\ldots , y_m)$ depend only on $\theta (x_0,\ldots ,x_k)$. Notice also that Theorem \ref{FVunital} is effective.
 
 Notice that equation (\ref{equivMain}) in Theorem \ref{FVunital} is  a statement only about the ring $R$. The left hand condition can be expressed as a condition on $R$ in terms of $B$-formulas. 
 
 It is worthwhile to compare a model $R$ of our axioms, and the full product $\prod_e R_e$   over all atoms $e$ of the Boolean algebra $\mathbb B$ of the idempotents of $R$,  which is also a model of the axioms. The obvious difference between $R$ and $\prod_e R_e$  is that $R$ may be countable, with $\mathbb B$ incomplete, whereas $\prod_e R_e$ will be uncountable when $\mathbb B$ is infinite and each $R_e$ has cardinality $\geq 2$. 
 
 
However, though both their Boolean algebras are atomic, and have the same atoms (i.e. there is a natural bijection between their atoms, and we pretend it is the identity), they will be different in general. The full product has Boolean algebra $\mathbb B_{\prod}$  the power set of the set of atoms, and it is complete.  The Boolean algebra $\mathbb B$ associated to $R$ need not be complete.

We construct a monomorphism from $\mathbb B$ to $\mathbb B_{\prod}$ as follows.  Let $r$ be a idempotent  of $R$. Then in either of the atomic Boolean algebras any nonzero element is the supremum of the atoms $e\leq r$. So we naturally match $r$ to the supremum in $B_{\prod}$ of the atoms $e\leq r$ (identifying the atoms of $\mathbb B$ with those of $\mathbb B_{\prod}$). This defines a map $F$ which is clearly an injective homomorphism of Boolean rings, and therefore an injective morphism of Boolean algebras (which we also call $F$.)
It is a well known fact (see e.g. \cite{GH}) that 
for atomic Boolean algebras $B_1$ and $B_2$, say with infinitely many atoms, any embedding is elementary if it preserves the notions of atom and union of at most $n\in \mathbb N$ atoms. 
This is certainly true in our situation, and so $F$ is 
elementary as a map from  $\mathbb B$ to the complete Boolean algebra  $\mathbb B_{\prod}$.

 Also the Boolean value of sentences is the same  in both algebras.  (Note that we have not tried yet to embed $R$ into the full product, this will be analyzed in future work).
At any rate, using Theorem \ref{FVunital}, this is enough to show that sentences    have the same truth value in the two rings, $R$ and $ \prod_{e \mbox{ {\tiny atom}}} R_e$, which is stated in the next corollary.

\begin{corollario}
\label{corFVunital}
Let $R$ as above. Then $$R\equiv \prod_{e \mbox{ {\tiny atom}}} R_e.$$
\end{corollario}

\proof
Notice that $R$ and the product of the $R_e$,  for $e$ a minimum idempotent, have the same Boolean algebra of idempotents, and the same stalks.  It is enough to apply the above theorem to sentences, and we use that $\mathbb B$ is an elementary substructure of $\mathbb B_{\prod}$. 
\endproof

\section{Exotic examples from Arithmetic}
\subsection{}
Since the 1930's one has studied nonstandard models of arithmetic, in several closely related formalisms. One formalism involves $+,\cdot, 0,1$ and $<$. The models for this formalism satisfy first-order Peano axioms, $PA$. See \cite{Kaye} for a clear presentation. The formalism classically applies to the semiring $\mathbb N$, and the model theory is pervaded by G\"odel Incompleteness. The other main formalism is inspired by the ring $\mathbb Z$, with the primitives above, as well as $-$. There is little difference logically (in terms of interpretability) between the two settings. In both the order $<$ is easily definable from the other relations via some simple axioms and Lagrange's Theorem. We have proved, using Ax's fundamental work \cite{Ax} of 1968 that although both set-ups are pervaded by G\"odel Incompleteness, they seem rather immune to any Incompletenss concerning the algebra of residue rings $\mathcal M/n\mathcal M$, where $n$ is a (possibly) nonstandard power of a prime in $\mathcal M$, for $\mathcal M$  a model of $PA$ (see \cite{PAresrings}). The analysis of the case $n$ composite is much harder, but can be handled by the earlier analysis in the present paper. The key is to formulate and prove a {\it factorization theorem} for $\mathcal M/n\mathcal M$ giving  a meaning to $$\mathcal M/n\mathcal M= \prod  \mathcal M/q\mathcal M$$ where $q$ varies among  the maximal prime powers dividing $n$. Using the refined analysis of the rings $\mathcal M/q\mathcal M$ in \cite{PAresrings} and 
the previous results in Section 2 we are able to understand the first order theory of the rings $\mathcal M/n\mathcal M$ for any $n\in \mathcal M$.

\subsection{}
We will work with the formalism for the semiring. If $\mathcal M$ is a semiring then it is easy to interpret the associated ring $\mathcal M^{\bullet}$ inside $\mathcal M$. Note that if $j$ is in the semiring $\mathcal M$ then $\mathcal M/j\mathcal M$ is naturally isomorphic to $\mathcal M^{\bullet}/j\mathcal M^{\bullet}$. Thus, for intelligibility we conflate  $\mathcal M$ and $\mathcal M^{\bullet}$ in what follows. 


G\"odel primitive recursive coding of syntax extends simply and uniformly to $\mathcal M$, giving a meaning to (not necessarily finite) formulas, and allowing certain inductive arguments about semantics. Details of such arguments and results are to be found in \cite{Kaye}. For us the crucial notions coded uniformly allowing us to prove properties about them in $PA$ are
\begin {enumerate}
\item
formulas, and specifically 
 $\Delta_0$-formulas;
 \item 
 satisfaction in $\mathcal M/n\mathcal M$, i.e. for any formula $\psi(v_0, \ldots, v_k)$ in $\mathcal L_{ring}$ and any $k+1$-tuple $c_0, \ldots, c_k$ in $\mathcal M$ we can code $$\mathcal M/n\mathcal M \models \psi(c_0+n\mathcal M, \ldots, c_k+n\mathcal M).$$ 
\end{enumerate}

\noindent 
Here $\mathcal M/n\mathcal M$ is represented as in elementary number theory as $[ 0, n-1]$ with $+$ and $\cdot$ defined via remainders modulo $n$.

The essential point is (2). The indespensable uniformity is that there is a single formula (not $\Delta_0$) $Sat(x,y,z)$, so that 
$Sat ( <\psi>,n, <\overline{c}> )$ expresses for all formulas $\psi$, and $ n,  \overline{c} \in \mathcal M$ that $$\mathcal M/n\mathcal M \models \psi(c_0+n\mathcal M, \ldots, c_k+n\mathcal M),$$  where $<\psi>$ is a code for the formula $\psi$, $n\in \mathcal M$, and $<\overline{c}>$ is a code for the $k+1$-tuple $\overline{c}=c_0, \ldots, c_k$ of elements of $\mathcal M$. This implies that $\mathcal M/n\mathcal M$ is recursively saturated. Note that not all rings elementarily equivalent to $\mathcal M/n\mathcal M$ are recursively saturated. 

Making this precise seems unecessary, as such codings and the notion of recursive saturation are ubiquitous in the subject.

We may interpret $\mathcal M/n\mathcal M$ inside  $\mathcal M$ in the usual way using the $\Delta_0$-relation of congruence mod $n$, picking  by $\Delta_0$-induction the smallest positive representative in each coset modulo the ideal $n\mathcal M$. We may also define on $[0,n-1]$ a $\Delta_0$-ring structure using Euclidean division, and this is what we then take as $\mathcal M/n\mathcal M$.  Thus we have a bounded $\Delta_0$-set and $\Delta_0$-relations $\oplus$, $\odot$, $\ominus$, with usual $0,1$ (notice that these functions are not related to those defined on the Boolean ring in Section 2.1). The aim is to show that $\mathcal M/n\mathcal M$ is a model of $AXIOMS$.  



Now we recall the factoriality properties of $\mathcal M$.  For fixed $n$ and prime $p$ in $\mathcal M$ we can define in a $\Delta_0$-way  the maximal power of $p$ dividing $n$, and we can prove that it is unique. The set $Q$ of maximal prime powers dividing $n$ is a bounded $\Delta_0$-set, and a suitable
version of the fundamental theorem of arithmetic can be proved in $PA$, giving a meaning to {\it $n$ is the product of prime powers}  in any non standard model of $PA$. Moreover, there is a clear $\Delta_0$-version of the Chinese Remainder Theorem ($\Delta_0$-CRT) which we state now for convenience since we will appeal to it in many proofs (see \cite{DA1} and  \cite{solPellLoc} for details and applications).
\begin{thm}
Let $\mathcal M$ be a model of $PA$ and $A$ a bounded $\Delta_0$-definable set in $\mathcal M$. Let $f$ and $r$ be $\Delta_0$-functions such that $f(a_1), f(a_2)$ are pairwise coprime for all $a_1, a_2 \in A$ and $a_1\not= a_2$, and
 $r(a) < f(a)$ for all $a \in A$. Suppose there exists $w\in \mathcal M$ divisible by all elements of $f(A)$. Then there exists $u <\prod_{a\in A} f(a)$ such that $u \equiv r(a) (\mbox{mod }f(a))$ for all $a \in A$.
\end{thm}
A much more general result holds giving a $\Delta_0$-definition of the product relation $y=\prod_{i\leq n} F(n)$, for $F$ a $\Delta_0$-definable function, see \cite{BDA}. This has been also used by the authors in the analysis of the multiplicative group of the residue field $\mathcal M/p\mathcal M$ where $\mathcal M$ is a model of a weak fragment of Peano Arithmetic, and $p$ is a prime in $\mathcal M$, see \cite{PAresfield}. 

An alternative way of considering  $\mathcal M/n\mathcal M$, which turns out to be convenient in our setting, is to identify $\mathcal M/n\mathcal M$ as a set of $\Delta_0$-definable functions. We now specify this. 

To each $x\in \mathcal M/n\mathcal M$ we assign the $\Delta_0$-function $\overline{x}$ defined on $Q$, the set of maximal prime powers $q$ dividing $n$, as follows 
$$
\begin{array}{ccccc}
\overline{x} & : & Q & \rightarrow & \mathcal M/n\mathcal M \\
&       &        q &         \mapsto &        \overline{x}(q)=rem(x,q)
\end{array}
$$
 where $rem(x,q)$ stands for the reminder of the division of $x$ by $q$ in $\mathcal M$. 
We now identify the idempotents of $\mathcal M/n\mathcal M$, i.e. the elements of the Boolean algebra $\mathbb B_n$ associated to $\mathcal M/n\mathcal M$.  Suppose $r\in \mathcal M/n\mathcal M$ is an idempotent then  $r^2-r\equiv 0(\mbox{mod }n)$ in $\mathcal M$. Therefore, for each $q$ maximal  prime power dividing $n$, $q$ divides $r^2-r$, so $q$ divides one, and only one, of $r$ or $r-1$. Conversely, if every $q$ in $Q$ divides either $r$ or $r-1$ then clearly $r$ is an idempotent in  $\mathcal M/n\mathcal M$. So,  unless $n$ is a power of a prime (in which case it is a local ring, and so connected) $\mathcal M/n\mathcal M$ is not connected. 

The minimal idempotents, i.e. the atoms in $\mathbb B_n$,  are exactly the $r$ such that for a unique $q\in Q$, $r\equiv 1($mod $q)$,  and $r\equiv 0($mod $q')$   for all other $q'\in Q$. By $\Delta_0$-CRT there are such $r$ for each $q$. Note that  $\mathcal M/n\mathcal M$ and the full product $\prod_{q\in Q}  \mathcal M/q\mathcal M$ have the same idempotents. 

We have then identified each idempotent in $\mathbb B_n$  with  a $\Delta_0$-definable  subset of $Q$, the set of maximal power of primes dividing $n$. The atoms  of $\mathbb B_n$  correspond to the subsets $Q-\{ q\}$, as $q$ varies in $Q$. To each atom  $e$ in $\mathbb B_n$  we associate the unique prime power $q_e$ such that $q_e\not|e$.  
We are aware that $q_e$ is not the best notation for expressing the 1-1- correspondence between atoms of $\mathbb B_n$ and prime powers dividing $n$. We hope that no confusion will arise.

 It is easy to show now that $(\mathcal M/n\mathcal M)_e$ is isomorphic $\mathcal M/q_e\mathcal M$. So, the elements of the localized ring at $e$ can be identified with elements of $\mathcal M$ which are $<q_e\leq n$.

Now we can verify that the ring $\mathcal M/n\mathcal M$ satisfies our axioms. 

\

\noindent
{\bf Axiom 1.} $\mathbb B_n$ is atomic. 

Let $r\in \mathbb B_n$, $r\not= 0$. There exists  $q\in Q$,  such that $r\equiv 1($mod $q)$, if not then $r=0$. Let $e$ be the atom such that  $e \equiv 1($mod $q)$ and $e \equiv 0($mod $q')$ for any $q'\in Q$ and $q'\not=q$. So, $re\equiv e(\mod n)$ which implies $e  \leq r$.  \hfill $\Box$
\endproof

\medskip

From now on an atom will always be an atom of $\mathbb B_n$. Let $e$ be an atom. By Lemma \ref{localized} the rings  $(\mathcal M/n\mathcal M)_e$ and  $ e(\mathcal M/n\mathcal M)$ are isomorphic via the map $\frac{r}{s}\mapsto er$ for $s=1,e$. Using this we can easily interpret  $(\mathcal M/n\mathcal M)_e$ inside $\mathcal M/n\mathcal M$ as the ring $ e(\mathcal M/n\mathcal M)$ which is definable (uniformly in $e$) inside  $\mathcal M/n\mathcal M$. We will also identify the ring $e(\mathcal M/n\mathcal M)$ with $\mathcal M/q_e\mathcal M$, where $q_e$ is the unique prime power dividing $n$ such that  $q_e\not | e$. 



\medskip
\noindent
{\bf Axiom 2.} Let $\theta (v_0, \ldots ,v_k)$ be any $\mathcal L_{rings}$-formula, and $f_0,\ldots, f_k\in \mathcal M/n\mathcal M$. We have to show that  there is a unique element $\llbracket   \theta (f_0, \ldots, f_n)   \rrbracket$ in $\mathbb B_n$ such that 
$$e(\mathcal M/n\mathcal M) \models  \theta (ef_{0}, \ldots ,ef_{k}) \mbox{ \hspace{.1in} iff  \hspace{.1in}} e \leq \llbracket   \theta (f_0, \ldots, f_n)   \rrbracket$$ for all atoms $e$ in $\mathbb B_n$. Let $$X=\{ e: e \mbox{ atom  and }  e(\mathcal M/n\mathcal M)\models  \theta (ef_{0}, \ldots ,ef_{k})\}.$$ Clearly, $X$ is a bounded $\Delta_0$-set of $\mathcal M$ since $\mathcal M/n\mathcal M$ lives in a $\Delta_0$-way on $[0,n-1]$, and so the relativization of  $\theta (f_0, \ldots ,f_k)$ gets coded in a $\Delta_0$-way. Now, $\Delta_0$-CRT implies the existence of the idempotent $\gamma $ in $\mathcal M/n\mathcal M$ such that  $\gamma\equiv 1($mod $q_e)$ for all $e\in X$,  
and  $\gamma\equiv 0($mod $q)$ for all other prime powers $q$ in $Q$. So, we define $\bigvee X=\gamma =\llbracket   \theta (f_0, \ldots, f_n)   \rrbracket$. \hfill $\Box$

\medskip
As we have just shown in our setting the join of the set $$\{ e: e \mbox{ is an atom and  }e(\mathcal M/n\mathcal M)\models  \theta (ef_{0}, \ldots ,ef_{k})\}$$ always exists, and  $$\bigvee \{ e: e \mbox{ is an atom and  }  e(\mathcal M/n\mathcal M)\models  \theta (ef_{0}, \ldots ,ef_{k})\}=\llbracket   \theta (f_0, \ldots, f_k)   \rrbracket.$$

\bigskip

\noindent
{\bf Axiom 3.}  Let $\exists v_{k+1}\theta (v_0, \ldots ,v_k, v_{k+1})$ be an $\mathcal L_{rings}$-formula, and $f_0,\ldots, f_k\in \mathcal M/n\mathcal M$. We need to show that  there is $g\in \mathcal M/n\mathcal M$ such that 
$$ \llbracket  \exists w \theta (\overline{f},w)     \rrbracket \leq \llbracket   \theta (\overline{f},g)     \rrbracket 
$$
where $\overline{f}=f_0, \ldots, f_n$. 

Let $$X=\{ e: e  \mbox { is an atom and }   e(\mathcal M/n\mathcal M)\models  \exists v_{k+1}\theta (ef_{0}, \ldots ,ef_{k},v_{k+1})\},$$ and $\bigvee X= \llbracket  \exists w \theta (\overline{f},w)     \rrbracket $ in $\mathcal M/n\mathcal M $.  Then by $\Delta_0$-induction for each $e$ in $X$ we take the minimum $h\in \mathcal M/n\mathcal M$ such that  $ e(\mathcal M/n\mathcal M)\models \theta (ef_{0}, \ldots ,ef_{k},eh)$. We will denote it by $h(e)$. 
Now we use again the $\Delta_0$-CRT in order to determine $g \in \mathcal M/n\mathcal M$ such that $g\equiv h(e)(\mbox{mod }q_e)$ for all atoms $e\in X$, and $g\equiv 0(\mbox{mod }q)$ for all other prime powers $q$ in $Q$. 

The set $$Y=\{ e:  e  \mbox { is an atom and }    e(\mathcal M/n\mathcal M)\models \theta (ef_{0}, \ldots ,ef_{k},eg)\}$$ is $\Delta_0$-definable and bounded, and arguing as before we can show that in $ \mathcal  M/n\mathcal M$ there exists the  supremum of $Y$ which we denote by $\llbracket   \theta (\overline{f},g)     \rrbracket$. From $X\subseteq Y$ it follows $$\llbracket  \exists v_{k+1} \theta (\overline{f},v_{k+1})     \rrbracket   \leq \llbracket   \theta (\overline{f},g)     \rrbracket.$$. \hfill $\Box$


\bigskip

\noindent
{\bf Axiom 4.} Let $\theta (v_0,\ldots , v_k)$ be an atomic $\mathcal L_{rings}$-formula, and $\overline{f}=f_1, \ldots , f_k\in \mathcal M/n\mathcal M$. We have to show that      $ \mathbb B_n \models  \llbracket   \theta (\overline{f})   \rrbracket=1 $ iff $\mathcal M/n\mathcal M \models  \theta (\overline{f}) $.

Now, $ \mathbb B_n \models  \llbracket   \theta (\overline{f})   \rrbracket=1 $ \   \   iff  \  \     $ \{ e: e\mbox{ atom  and } e(\mathcal M/n\mathcal M)\models \theta (e\overline{f})\} $ coincides with the set of atoms of $\mathcal M/n\mathcal M$   \  \  iff  \   \
$\mathcal M/q_e\mathcal M\models \theta (\overline{f}_{q_e}) $  for all atoms $e$ \   \   iff  \  \    $\mathcal M/n\mathcal M \models  \theta (\overline{f}) $. Notice that the last equivalence follows from the fact that for a tuple to satisfy an atomic formula corresponds to be a zero of a polynomial.

\bigskip

\noindent
{\bf Axiom 5.} The proof of this axiom in $\mathcal M/n\mathcal M$ follows the line of  the proof of  Lemma \ref{axiom5}.  The only significant modification is in the proof of 2) $\Rightarrow$ 1) where the element $g$ is obtained using $\Delta_0$-CRT  without appealing to the axiom of choice. 

\smallskip

Now, arguing as in the proof of Corollary \ref{corFVunital} we have 

\begin{thm}
\label{main}
$\mathcal M/n\mathcal M\equiv \prod_{q\in Q}  \mathcal M/q\mathcal M$.
\end{thm}

\medskip

 In a forthcoming paper  we will use Theorem \ref{main} to describe thoroughly the model theory of the rings $\mathcal M/n\mathcal M$ for any $n\in \mathcal M$.

\medskip

Derakhshan and Macintyre in \cite{DM-ad2} adapt our methods to the ad\`elic situation, at the cost of modifying the preceeding axioms. The main difference is that the Boolean algebras of idempotents are enriched with extra predicates identifying those elements which are join of a finite sets of atoms.

\medskip
\noindent
{\bf Acknowledgements.} The authors are very grateful to MSRI in Berkeley where this paper was completed during the program Definability, Decidability and Computability in Number Theory in July-August 2022. The first author was partially supported  by PRIN 2017-Mathematical Logic: Models, Sets, and Computability. The second author was supported by a Leverhulme Emeritus Fellowship. 

We are grateful to Jamshid Derakhshan for many stimulating and helpful discussions.

\bibliographystyle{plain}
\bibliography{biblioU}

\def\Dbar{\leavevmode\lower.6ex\hbox to 0pt{\hskip-.23ex \accent"16\hss}D}
\begin{thebibliography}{10}

\bibitem{Ax}
James Ax.
\newblock The elementary theory of finite fields.
\newblock {\em Ann. of Math. (2)}, 88:239--271, 1968.

\bibitem{Bell}
J.~L. Bell.
\newblock {\em Boolean-valued models and independence proofs in set theory},
  volume~12 of {\em Oxford Logic Guides}.
\newblock The Clarendon Press, Oxford University Press, New York, second
  edition, 1985.

\bibitem{BDA}
Alessandro Berarducci and Paola D'Aquino.
\newblock {$\Delta_0$}-complexity of the relation {$y=\prod_{i\leq n}F(i)$}.
\newblock volume~75, pages 49--56. 1995.
\newblock Proof theory, provability logic, and computation (Berne, 1994).

\bibitem{PAresrings}
P.~D'Aquino and A.~Macintyre.
\newblock The model theory of residue rings of models of {P}eano {A}rithmetic:
  the prime power case.
\newblock {\em arXiv:3582301}, 2020.

\bibitem{DA1}
Paola D'Aquino.
\newblock Local behaviour of the {C}hebyshev theorem in models of {${\rm
  I}\Delta_0$}.
\newblock {\em J. Symbolic Logic}, 57(1):12--27, 1992.

\bibitem{solPellLoc}
Paola D'Aquino.
\newblock Solving {P}ell equations locally in models of {${\rm I}\Delta_0$}.
\newblock {\em J. Symbolic Logic}, 63(2):402--410, 1998.

\bibitem{PAresfield}
Paola D'Aquino and Angus Macintyre.
\newblock Non-standard finite fields over {$I\Delta_0+\Omega_1$}.
\newblock {\em Israel J. Math.}, 117:311--333, 2000.

\bibitem{DM-ad2}
J.~Derakhshan and A.~Macintyre.
\newblock Decidability problems for adele rings and related restricted
  products.
\newblock {\em arXiv:1910.14471}, 2019.

\bibitem{FV1959}
S.~Feferman and R.~L. Vaught.
\newblock The first order properties of products of algebraic systems.
\newblock {\em Fund. Math.}, 47:57--103, 1959.

\bibitem{GH}
Steven Givant and Paul Halmos.
\newblock {\em Introduction to {B}oolean algebras}.
\newblock Undergraduate Texts in Mathematics. Springer, New York, 2009.

\bibitem{Kaye}
Richard Kaye.
\newblock {\em Models of {P}eano arithmetic}, volume~15 of {\em Oxford Logic
  Guides}.
\newblock The Clarendon Press, Oxford University Press, New York, 1991.
\newblock Oxford Science Publications.

\bibitem{MacResField}
A.~Macintyre.
\newblock Residue fields of models of {${\rm P}$}.
\newblock In {\em Logic, methodology and philosophy of science, {VI}
  ({H}annover, 1979)}, volume 104 of {\em Stud. Logic Found. Math.}, pages
  193--206. North-Holland, Amsterdam, 1982.

\bibitem{Mos}
A.~Mostowski.
\newblock On direct products of theories.
\newblock {\em J. Symbolic Logic}, 17:1--31, 1952.

\bibitem{Weispf}
V.~Weispfenning.
\newblock Quantifier elimination and decision procedures for valued fields.
\newblock In {\em Models and sets ({A}achen, 1983)}, volume 1103 of {\em
  Lecture Notes in Math.}, pages 419--472. Springer, Berlin, 1984.

\bibitem{BorisComm}
B.~Zilber.
\newblock The semantics of the canonical commutation relation.
\newblock {\em arXiv:1604.07745}, 2016.

\end{thebibliography}
\end{document}